\documentclass[12pt]{article}
\usepackage{amssymb,latexsym,amsmath}
\def\C{\mathbb C}
\def\R{\mathbb R}
\def\N{\mathbb N}
\def\Z{\mathbb Z}

\newtheorem{thm}{Theorem}[section]

\newtheorem{lem}{Lemma}[section]

\newtheorem{definition}{Definition}[section]

\begin{document}
\sffamily

\title{The  Schwarzian derivative\\ and the Wiman-Valiron property}
\author{J.K. Langley}
\maketitle

\begin{abstract}
Let $f$ be a transcendental meromorphic function in the plane such that
$f$ has finitely many critical values, the
multiple points of $f$ have bounded multiplicities, and
the inverse function of $f$ has finitely many transcendental singularities.
Using the Wiman-Valiron method it is shown that if the Schwarzian derivative $S_f$ of $f$ is transcendental then
$f$ has infinitely many multiple points,
the inverse function of $S_f$ does not have a direct transcendental singularity over $\infty$, and
$\infty$ is not a Borel exceptional value of $S_f$. The first of these conclusions was proved by Nevanlinna and Elfving
via a fundamentally different method.\\
A.M.S. MSC 2000: 30D35.
\end{abstract}

\section{Introduction}\label{intro}

The singularities of the inverse of a transcendental meromorphic function $f$ in the plane
play a pivotal role in function theory and complex dynamics,
and were classified by Iversen \cite{Erem,Iversen} as follows. An algebraic singularity of $f^{-1}$ 
over a value $a$ in the extended complex plane
$\C^\infty$ arises from a multiple point $z \in \C$ of
$f$ with $f(z) = a$, that is, a multiple pole of
$f$ or a zero of $f'$, while a path $\gamma$ tending to infinity on which $f(z)$ tends to  $a \in \C^\infty$ 
determines a transcendental singularity of $f^{-1}$ over $a$. Paths $\gamma_1, \gamma_2$ with the same asymptotic value $a$ 
determine distinct transcendental singularities over $a$ if there exist $t > 0$ and disjoint components
$C_1$, $C_2$ of the set $\{ z \in \C : \chi( f(z), a) < t \}$ with $\gamma_j \setminus C_j$ bounded: here 
$\chi(a, b)$ denotes the chordal metric on $\C^\infty$. 
The starting point of the present paper 
is the following theorem, the restriction to the plane meromorphic setting of results of Elfving \cite{elf} and Rolf Nevanlinna \cite{Nev2,Nev}.

\begin{thm}
[\cite{elf,Nev2,Nev}]
\label{thmA}
Let $f$ be a transcendental meromorphic function in the plane satisfying both of the following conditions.\\
(i) The function
$f$ has finitely many multiple points.\\
(ii) The inverse function $f^{-1}$ 
has finitely many transcendental singularities. \\
Then the Schwarzian derivative
\begin{equation}
 \label{Sfdef}
S_f(z) = \frac{f'''(z)}{f'(z)} - \frac32 \left( \frac{f''(z)}{f'(z)} \right)^2 
\end{equation}
is a rational function. 
\end{thm}
In fact, the papers \cite{elf,Nev2} consider simply connected
Riemann surfaces spread over the sphere with finitely many algebraic and logarithmic
singularities (see ~\cite{Er2}), and together show that all such surfaces 
are uniformised by meromorphic functions in the plane with rational Schwarzian
\cite[p.57]{elf}, and in particular are parabolic. On the other hand,  for certain choices of  angles in the construction of the Schwarz triangle
functions,  the inverse has infinitely many singularities, each of them
algebraic and lying over one of three values, and  
the Riemann surface is hyperbolic \cite[pp.298-9]{Nev} (see also \cite[Ch. 7]{Nehari2}). 

The converse of Theorem \ref{thmA} is also true, in the sense that if $f$ is a transcendental  meromorphic function in the plane and 
$S_f$ is a rational function then $f$ satisfies (i) and (ii). Moreover, the Schwarzian $S_f$ is intimately linked
to differential equations since $f$ may be written locally as a quotient of solutions of the equation
\begin{equation}
 \label{1}
w'' + A(z) w = 0,
\end{equation}
where $2A = S_f$. In the case where $A \not \equiv 0$ is a polynomial of degree $n$,  asymptotic integration
\cite{Hil1,Hil2} shows that the quotient $f$ of linearly independent solutions of 
(\ref{1}) is a meromorphic function in the plane, of order of growth $(n+2)/2$, with no multiple points and $n+2$ transcendental
singularities of the inverse \cite{FNev}. These functions $f$ are  extremal for Nevanlinna's deficiency relation
$\sum \delta(a, f) \leq 2$, while functions of finite lower order
which have maximal deficiency sum, or sufficiently few multiple points, have asymptotic behaviour which closely resembles that of the quotient of two
solutions of (\ref{1}) \cite{Dra1,Er2}. For the definitions of the order of growth and the Nevanlinna deficiency see
\cite{Hay2, Nev}. 

The proof of Theorem \ref{thmA} in \cite{elf,Nev2,Nev} uses the method of line complexes \cite{GO,Nev} to approximate $S_f$ by the Schwarzian derivative
of a rational function, and the present paper will focus on 
the following alternative
approach. Suppose that $f$ is a transcendental meromorphic function in the plane with finitely many multiple points, 
but that $S_f$ is transcendental. Then $S_f$ has finitely many poles, so the Wiman-Valiron theory \cite{Hay5} gives small
neighbourhoods on which $A(z) = S_f(z)/2$ behaves like a power of $z$ (see \S\ref{WV}), and local asymptotic representations for solutions of
(\ref{1}) on these neighbourhoods were developed in \cite{La5}. Since the powers arising from the Wiman-Valiron theory tend to
infinity, it is reasonable to anticipate that these local asymptotics are incompatible with the existence of only finitely many 
transcendental singularities of $f^{-1}$. It turns out that this approach does indeed work and, although  it
may not be simpler than the proof of Theorem \ref{thmA} in \cite{elf,Nev2,Nev}, it does deliver some additional conclusions 
as follows. 

\begin{thm}\label{thm1}
Let $M$ be a positive integer and let $f$ be a transcendental meromorphic function in the plane
with Schwarzian derivative $S_f =2A$ and the following two properties:\\
(i) $f$ has finitely many critical values and
 all multiple points of $f$ have multiplicity at most $M $;\\
(ii) the inverse function of $f$ has finitely many transcendental singularities.\\
If $A$ is transcendental then the following three conclusions hold.\\
(a) The function $f$ has infinitely many multiple points.\\
(b) The inverse function of $A$ does not have a direct transcendental singularity over $\infty$.\\
(c) The value $\infty$ is not Borel exceptional for $A$, that is, 
\begin{equation*}
 \label{conc1}
\limsup_{r \to \infty} \frac{ \log N(r, A)}{\log r} =  \limsup_{r \to \infty} \frac{ \log T(r, A)}{\log r} .
\end{equation*}
 \end{thm}

Here \cite{BE,Nev} a transcendental singularity, arising from a path $\gamma$ tending to infinity 
on which a transcendental meromorphic function $h$ 
tends to 
$a \in \C^\infty $, is said to be direct if there exist $t > 0$ and a 
component
$C$ of the set $\{ z \in \C : \chi( h(z), a) < t \}$ such that $h(z) \neq a$ on $C$ and  $\gamma \setminus C$ is bounded. 
By Iversen's theorem the inverse of a transcendental entire function always has a direct transcendental singularity over $\infty$,
but so has, for example, that of $e^{-z^2} \cot z$.

Note that in Theorem \ref{thm1} it is not assumed initially that $f$ has finitely many 
multiple points, but because of (i) and (ii) the inverse function has finitely many singular values, so that $f$ 
belongs to the Speiser class $S$ \cite{EL}. 
Conclusion (a) is not new, and merely replicates that of Theorem \ref{thmA}, albeit with a different proof,
but 
(b) and (c) are new as far as the author is aware. 
There are abundant examples of functions satisfying (i) and (ii) with transcendental Schwarzian, such as $\cos z$ or the Weierstrass 
doubly periodic function $\wp(z)$; one could also take $\wp \circ g$ for any entire $g$
with no finite asymptotic values such that the critical points of $g$ have bounded multiplicities
and all but finitely many of them are mapped by $g$ to multiple points of $\wp$.

\section{The Wiman-Valiron property}\label{WV}


\begin{definition}
 \label{wvprop}
Let $R \geq 0$ and let the function $A$ be meromorphic in  $R \leq |z| < \infty$ with an essential singularity at infinity.
Then $A$ has the Wiman-Valiron property if there exist an unbounded set $F_1 \subseteq [R, \infty)$ 
and functions $N: F_1 \to [2, \infty )$ and $\phi: F_1 \to [1, \infty )$ 
which satisfy the following two conditions.\\
(i) The functions $N$ and $\phi$ have
\begin{equation}
 \label{Nlim}
\lim_{r \to \infty, r \in F_1 } N(r) = \infty, \quad
\lim_{r \to \infty, r \in F_1} \frac{\phi(r)}{N(r)} = 0 ,
\quad \lim_{r \to \infty, r \in F_1} \frac{\phi(r)}{\log N(r)} = \infty .
\end{equation}
(ii) For each $r \in F_1$ there exists $z_r$ with
\begin{equation}
 \label{zrdef}
|z_r| = r, \quad \phi(r) = o( \log^+ |A(z_r)| ) ,
\end{equation}
such that
\begin{equation}
 \label{Nlim2}
A(z) \sim \left( \frac{z}{z_r} \right)^{N(r)} A(z_r) 
\end{equation}
on  $D(z_r, 8)$, uniformly as $r \to \infty$ in $F_1$, where  
\begin{equation}
 \label{Nlim1}
D(z_r, L) = 
\left\{  z_r e^\tau : \, 
| {\rm Re} \, \tau  | \leq L \frac{\phi(r) }{N(r)} \, , \,  
| {\rm Im } \, \tau | \leq L  \frac{\phi(r) }{N(r)} \, \right\}.
\end{equation} 
\end{definition}
Note that 
$$\log |A(z)| \geq \log |A(z_r)| - 8 \phi (r) - o(1) \sim \log |A(z_r)| $$
on $D(z_r, 8)$, by (\ref{Nlim}), (\ref{zrdef}) and (\ref{Nlim2}). 

The Wiman-Valiron property  is well known for  transcendental entire  functions $A$ \cite{Hay5}, in which case 
one may take $\phi(r) = N(r)^{1-\gamma}$ for any $\gamma $ with $1/2 < \gamma < 1$, and $z_r$ with
$|z_r| = r$ and $|A(z_r)| \sim M(r, A)$; here $M(r, A)$ is the maximum modulus,
$N(r)$ is the central index of $A$ and 
$[1, \infty ) \setminus F_1$ has finite logarithmic measure. With very minor modifications, 
the same estimates hold when $A$ is transcendental meromorphic with finitely many poles, which is automatically the case by
(\ref{Sfdef}) if $2A = S_f$ is transcendental and $f$ is meromorphic with finitely many multiple points. It will be 
noted in Lemma \ref{lemwv1} that any transcendental meromorphic function in the plane for which $\infty$ is a Borel 
exceptional value has the Wiman-Valiron property. More significantly, it was proved in \cite{BRS}
(see in particular \cite[Theorems 2.1 and 2.2]{BRS} and \cite[Lemma 6.10]{BRS}) that an almost exact counterpart of
the Wiman-Valiron theory holds for any transcendental meromorphic function $A$ in the plane for which the inverse function has a direct transcendental
singularity over 
$\infty$. Thus Theorem \ref{thm1} follows from the next result.

\begin{thm}
 \label{thm2}
Let $M$ be a positive integer and let $f$ be a transcendental meromorphic function in the plane with 
Schwarzian derivative $S_f =2A$ and the following two properties:\\
(i) $f$ has finitely many critical values and
all multiple points of $f$ have multiplicity at most $M $;\\
(ii) the inverse function of $f$ has finitely many transcendental singularities.
\\
If $A$ is transcendental then $A$ does not have the Wiman-Valiron property. 
\end{thm}

\section{Local asymptotic integration}

Let the function $A$ be meromorphic in  $R \leq |z| < \infty$ with an essential singularity at infinity, and assume that $A$ 
has the Wiman-Valiron 
property. Take $r \in F_1$ and $z_r$ as in Definition \ref{wvprop}, with $r$ large, and set 
\begin{equation*}
N = N(r), \quad 
w_r = z_r \exp \left(  \frac{ -4 \phi(r)}N \right).
\label{b6}
\end{equation*}
Then  (\ref{Nlim2}) may be written in the form
\begin{equation}
 \label{Nlim4}
A(z) =  \left( \frac{z}{w_r} \right)^{N} A(w_r) ( 1 + \mu (z) )^2 , 
\quad \mu (z) = o(1), 
\end{equation}
for $z$ in $D(z_r, 8)$, so that (\ref{Nlim}), (\ref{Nlim1}) and Cauchy's estimate for derivatives yield 
\begin{equation}
 \label{Nlim3}
\mu'(z) = o \left( \frac{N}{r\phi(r)}  \right), \quad 
\mu''(z) = o \left( \frac{N^2}{r^2\phi(r)^2 } \right) , \quad 
\frac{A'(z)}{A(z)} \sim  \frac{N}{z} \quad \hbox{and} \quad
\frac{A''(z)}{A(z)} \sim  \frac{N^2}{z^2}
\end{equation}
for $z$ in $D(z_r, 4)$. Again for $z$ in $D(z_r, 4)$, set 
\begin{equation}
 \label{Nlim5}
Z = \frac{2 w_r A(w_r)^{1/2}  }{N+2} + \int_{w_r}^z A(t)^{1/2} \, dt =
\frac{2 w_r A(w_r)^{1/2}  }{N+2} + \int_{w_r}^z  \left( \frac{t}{w_r} \right)^{N/2} A(w_r)^{1/2} ( 1 + \mu (t) ) \, dt .
\end{equation}
Now let $z \in D(z_r, 4)$ and let $\sigma_z$ be the path 
from $w_r$ to $z$ consisting of the radial segment from $w_r$ to
$z^* = w_r |z/w_r|$ followed by the shorter circular arc from $z^*$ to $z$, so that $\sigma_z$ has 
length $O(r \phi(r)/N  )$. Now $|w_r| \leq |t| \leq |z|$ on $\sigma_z$, and so (\ref{Nlim3})  and integration by parts along $\sigma_z$ yield
\begin{eqnarray*}
\int_{w_r}^z t^{N/2} \mu (t) \, dt &=&
o\left(  \frac{|z|^{(N+2)/2} }{N+2} \right) - \int_{w_r}^z o \left( \frac{N}{r\phi(r)}  \right) \, \frac{2  t^{(N+2)/2} }{N+2} \, 
  dt = o\left(  \frac{|z|^{(N+2)/2} }{N+2} \right)  .
\end{eqnarray*}
Hence, for $z \in D(z_r, 4)$, using (\ref{Nlim4}) and (\ref{Nlim5}), 
\begin{eqnarray}
\label{b7}
Z =Z(z) &\sim&  \frac{ 2 z^{(N+2)/2}A(w_r)^{1/2} }{ w_r^{N/2} (N+2) }  
\sim   \frac{2z A(z)^{1/2}  }{N+2}  \sim Z(z_r) \left( \frac{z}{z_r} \right)^{(N+2)/2} ,   \nonumber \\
\log \frac{Z(z)}{Z(z_r)} &=& \frac{N+2}{2} \log \frac{z}{z_r} + o(1) , \quad 
\frac{d \log Z}{d \log z} = \frac{z A(z)^{1/2}}{Z} \sim \frac{N+2}2 . 
\end{eqnarray}
Here the logarithms are determined so that $\log (z/z_r)$ vanishes at $z_r$. 
The following lemma is an immediate consequence of (\ref{b7}). 
\begin{lem}
 \label{wvlema}
Let $Q$ be a positive integer and let $r \in F_1$ be  large. 
Then $\log Z$ is a univalent function of $\log z$ on $D(z_r, 7/2)$ and there exist at least $Q$ 
pairwise disjoint
simple islands $H_q$ in $D(z_r, 3)$ mapped univalently by $Z$ onto
the closed logarithmic rectangle
\begin{equation}
J_1 = \left\{ Z: \, \left| \log \left|\frac{Z}{Z_r}\right| \right| \leq \phi(r) , \, 
| \arg Z | \leq \frac{\pi}4 \right\} , \quad
Z_r = \left| \frac{2 z_r A(z_r)^{1/2}  }{N+2} \right| \sim  | Z(z_r)|.
\label{b8}
\end{equation}
\end{lem}
\hfill$\Box$
\vspace{.1in}

The next step is to apply a local analogue of Hille's asymptotic method \cite{Hil1,Hil2}
developed in~\cite{La5}.
For $z$ in $H_q$ (where $1 \leq q \leq Q$)  and $Z = Z(z)$ in $J_1$, a solution $w(z)$
of (\ref{1}) is transformed via  
\begin{equation}
W(Z) = A(z)^{1/4} w(z) , \quad 
\frac{d^2W}{dZ^2} + (1 - F_0(Z))W = 0,
\quad
F_0(Z) =
\frac{A''(z)}{4A(z)^2} - \frac{5A'(z)^2}{16A(z)^3 } ,
\label{b10}
\end{equation}
and
$|F_0(Z)| \leq 3 |Z|^{-2}$ in $J_1$, using (\ref{Nlim3}) and (\ref{b7}). Since
$$
\log (Z_r \exp( - \phi(r) )) \geq 
\frac12 \log |A(z_r) | + \log r - \log N(r) - \phi (r) - O(1),
$$
(\ref{Nlim}) and (\ref{zrdef}) imply that
\begin{equation}
 \label{limZr}
\lim_{r \to \infty, r \in F_1 } \frac{Z_r}{ N(r)} = 
\lim_{r \to \infty, r \in F_1 }  Z_r \exp ( - \phi(r) ) = \infty .
\end{equation}
Hence \cite[Lemma 1]{La5} gives 
solutions $U_1(Z), U_2(Z)$ in $J_1$ of the differential  equation for $W$ appearing
in (\ref{b10}), as well as solutions $u_1(z), u_2(z)$ of (\ref{1}) in $H_q$, these 
satisfying
\begin{equation}
U_j(Z) \sim \exp ( (-1)^j i Z ), \quad u_j(z) = A(z)^{-1/4} U_j(Z) .
\label{b11}
\end{equation}
\hfill$\Box$
\vspace{.1in}

\section{The class $B$ and the Schwarzian derivative}

Throughout this paper,  $B(a, r)$ will denote the open Euclidean disc of centre $a$
and radius $r$, and $S(a, r)$  the corresponding
boundary circle. The chordal (spherical) metric
on the extended complex plane $\C^\infty$ will be written as
$$
\chi (z, w) =  \frac{| z - w |}{
\sqrt{  ( 1 + |z|^2 ) ( 1 + |w|^2 ) }} , \quad 
\chi (z, \infty ) =  \frac{1}{
\sqrt{   1 + |z|^2  }} .
$$
Suppose next that  $h$ is a transcendental meromorphic function with no asymptotic or
critical values in $0 < R < |w| < \infty $.
Then by a well known classification result \cite[p.287]{Nev}, every
component $C_0$ of the set
$\{ z \in \C : |h(z)| > R \}$ is simply connected, and there are
two possibilities. Either (i) $C_0$ contains 
one pole $z_0$ of $h$, of multiplicity $k$ say, in which case 
$h^{-1/k} $ maps $C_0$ univalently onto $B(0, R^{-1/k})$, or 
(ii) $C_0$ contains no pole of $h$, but instead a path
tending to infinity on which $h$ tends to infinity. In case (ii) the function $w = \log h(z)$ maps 
$C_0$ univalently onto the  half-plane ${\rm Re} \, w > \log R$,  
the singularity of $h^{-1}$ over $\infty$ is called logarithmic and, following \cite{BRS}, the component $C_0$ 
will be called a logarithmic tract. Logarithmic singularities over finite values are defined analogously and, if $h$ belongs to
the Speiser class $S$, all transcendental singularities of $h^{-1}$ are automatically logarithmic. 

The following lemma will be proved for functions in the Eremenko-Lyubich class $B$, that is, transcendental meromorphic functions
in the plane for which the set of finite singular values of the inverse function is bounded, so that any transcendental 
singularity over $\infty$ is again logarithmic.

\begin{lem}
 \label{slem1}
Let $1 \leq R < \infty$ 
and let $f$ be a transcendental meromorphic function in the plane such that all  finite critical and asymptotic values $w$ of $f$ satisfy 
$|w| \leq  R $. 
Let $U$ be a component of the set $\{ z \in \C :  |f(z)|  > R \}$ which contains no poles of $f$, and let $z_0 \in \C \setminus U$.
Then
\begin{eqnarray*}
 \label{slemest}
\left( \log \frac{|f(z)|}R \right)^2 &\leq&  T(z) = 32 \pi^2  | (z-z_0)^2 S_f(z)| + 16 \pi^2 + 12, 
\nonumber \\
\quad \chi( f(z) , \infty) &\geq& \frac1{2R} \exp\left( - \sqrt{ |T(z)| } \right)  \geq  \frac1{R}
\exp\left( - 20 |z-z_0| \sqrt{ | S_f(z)| } - 20 \right) , 
\end{eqnarray*}
for $z \in U$, where $ S_f(z)$ denotes the Schwarzian derivative of $f$.
\end{lem}
\textit{Proof.}
Write
$$
u = \psi(w) = \log ( f^{-1}( e^w) - z_0) , \quad e^w = f(e^u+z_0) = f(e^{\psi(w)} +z_0) , \quad z= e^u + z_0, 
$$
on $H = \{ w \in \C : {\rm Re} \, w > \log R \}$. Since $\infty \not \in f(U)$, the function $\psi$ is univalent on $H$.
Thus Koebe's quarter theorem \cite[p.3]{Hay9} and Nehari's univalence criterion \cite{Nehari} give 
\begin{equation}
 \label{slemest1}
|\psi'(w)| \leq \frac{4 \pi }{{\rm Re} \, w - \log R }  = \frac{4 \pi }{\log \frac{|f(z)|}R  } , \quad 
| S_\psi(w) | \leq \frac{6}{({\rm Re} \, w - \log R)^2 } =  \frac{6}{\left(\log \frac{|f(z)|}R \right)^2 } 
\end{equation}
on $H$. 
Next, the composition formula 
$$
S_{p \circ q} = S_p(q) (q')^2 + S_q
$$
for the Schwarzian and the fact that $S_p \equiv 0$ when $p$ is a M\"obius transformation give
\begin{eqnarray}
- \frac12  &=& S_{\exp} (w) = S_f( e^{\psi(w)}+z_0 ) e^{2 \psi(w)} \psi'(w)^2 + S_{\exp \circ \psi  } (w) \nonumber \\
&=& S_f(z) (z-z_0)^2 \psi'(w)^2  + S_{\exp} ( \psi(w) )  \psi'(w)^2 + S_\psi(w) \nonumber \\
&=& S_f( z)(z-z_0)^2 \psi'(w)^2 - \frac12 \psi'(w)^2 + S_\psi(w) .
\label{slemest2}
\end{eqnarray}
Now combining (\ref{slemest1}) and (\ref{slemest2}) with the inequality $\sqrt{A^2+B^2} \leq A+B$ for positive $A$ and $B$
completes the proof.
\hfill$\Box$
\vspace{.1in}

\section{Further estimates for functions in the class $S$}

Throughout this section let $f$ be a transcendental meromorphic function in the plane such that
all singular values of $f^{-1}$ lie in the finite set $\{ a_1, \ldots, a_\kappa \} \subseteq \C^\infty $. 
If $E$ is a subset of a simply connected subdomain of $\C \setminus \{ 0 \}$, define the 
\textit{logarithmic diameter} of $E$ to be 
$$
\sup_{u, v   \in E} \left|  \log u - \log v \right| . 
$$

\begin{lem}
 \label{component}
There exists a positive real number $R $, depending on $f$,  with the following property.  
Let $ \rho $  be small and positive and $\alpha \in \C^\infty$ and let $E \subseteq \{ z \in \C : R < |z| < \infty \}$ be a connected set 
such that $\chi( f(z), \alpha )  < \rho^2 $ for all $z \in E$. 
Then at least one of the following two conclusions holds:\\
(i) the set $E$ lies in a component of the set $\{ z \in \C : \chi( f(z), a_j) < 2 \rho \}$ for some $j$;\\
(ii) the set $E$ lies in a simply connected subdomain of $\C \setminus \{ 0 \}$ and 
has logarithmic diameter at most $400 \pi \rho $. 
\end{lem}
\textit{Proof.} Choose $x_0, x_1 \in \C$ and $R  $ such that 
$$
|f(x_0)| \leq \frac12, \quad |f(x_1)| \geq 2, \quad R \geq 2 ( |x_0| + |x_1| ). 
$$
With this choice of $R$, let $E$ be as in the statement of the lemma, assume that (i) does not hold, and choose $z_0 \in E$. 
Then  $w_0 = f(z_0)$ satisfies 
$\chi( f(z) , w_0) < 2 \rho^2 <  \rho $ for all $z \in E$, and 
$\chi( w_0, a_j) \geq   \rho $ for all~$j$. It may be assumed that $|w_0| \leq 1$. Then the branch 
of $g=f^{-1}$ mapping $w_0$ to $z_0$ extends to be analytic and univalent on the Euclidean disc
$$F = B(w_0, \rho) \subseteq  \{ w \in \C: \chi(w, w_0) <   \rho \} \subseteq B(0, 2) .$$
Now let $|w-w_0| < 10 \rho^2 < \rho/2$. Since  $\psi (w) = \log ( g(w) - x_1 )$ is  analytic and
univalent on $F$,  Koebe's quarter theorem gives
$$
|\psi'(w)| \leq \frac{8 \pi }{\rho}, \quad | \psi (w) - \psi (w_0) | \leq 80 \pi \rho . 
$$
Finally, let $z \in E$ and $w = f(z)$. Then $\chi(w, w_0) < 2 \rho^2$ and $|w-w_0| < 10 \rho^2 $ and so 
$$
z - x_1 = (z_0 - x_1) \exp( \psi(w) - \psi(w_0)) = (z_0-x_1)(1 + \delta (z)), \quad | \delta (z)| \leq 90 \pi \rho ,
$$
which leads easily to (ii). 
\hfill$\Box$
\vspace{.1in}

The next lemma requires additional assumptions on the function $f \in S$.

\begin{lem}
 \label{lemp2}
Assume that
all singular values $a_j$ of $f^{-1}$ are finite, with $|a_j - a_{j'}| \geq 2 \sigma > 0$ for $a_j \neq a_{j'}$. 
Assume further that all zeros of $f'$ have multiplicity at most $M-1$, for some fixed $M \in \N$.
Let $\rho$ be a positive real number, small compared to $\sigma$, and let $z_0$ be large with $| f(z_0)- a_j | < \rho $ for some $j$.
Assume that $a_j \in f(F)$,  where $F$ is that component of $G = \{ z \in \C : | f(z) -a_j  |  < \sigma  \}$
which contains $z_0$. 
Then $z_0$ lies in a simply connected component $D \subseteq \C \setminus \{ 0 \}$ 
of the set $\{ z \in \C : | f(z) - a_j | < \rho  \}$, and 
$D$ has logarithmic diameter at most 
$$
16 \pi \left( \frac{\rho}{\sigma}\right)^{1/M} .
$$
\end{lem}
\textit{Proof.} Assume without loss of generality that $a_j = 0$.
Then any component of $G$ which contains a zero of $f$ must be bounded. 
Since $0 \in f(F)$ and $z_0$ is large it may therefore be assumed that $0 \not \in F$. Now $f$ can be written in the form 
$f = g^m$ on $F$,
where $g: F \to B(0, \sigma^{1/m} )$ is conformal, and $m$ is a positive integer with $m \leq M$. 
Let $G: B(0, \sigma^{1/m} ) \to F$ be the inverse map of $g$. Then 
$$z_0 \in \{ z \in F : | f(z) | < \rho  \} = D = G(B(0, \rho^{1/m})).$$
Since $\psi(w) = \log G(w)$ is analytic and univalent
on $B(0, \sigma^{1/m})$, and $\rho/ \sigma$ is small, Koebe's quarter theorem yields, 
for $|w| \leq \rho^{1/m}$,
$$
|\psi'(w)| \leq \frac{8 \pi}{\sigma^{1/m}} , \quad \left| \log G(w) - \log G(0) \right| =
| \psi(w) - \psi(0) |    
\leq 8 \pi \left( \frac{\rho}{\sigma}\right)^{1/m} 
\leq 8 \pi \left( \frac{\rho}{\sigma}\right)^{1/M} .
$$
\hfill$\Box$
\vspace{.1in}

\section{Proof of Theorem \ref{thm2}}

Let $f$ be as in the hypotheses and assume without loss of generality that all singular values of $f^{-1}$ 
are finite but have modulus at least $2$. Take a small positive real number $\tau$. The 
finitely many transcendental singularities of $f^{-1}$
give rise to a (possibly empty) finite set of logarithmic tracts $W_k$, 
each a component of the set $\{ z \in \C : |f(z)-a_j| < \tau \}$ for some $j = j(k)$, such that 
$$
w= \log \frac1{f(z)-a_j} 
$$
maps $W_k$ conformally onto the half-plane ${\rm Re} \, w > - \log \tau $. By reducing $\tau$ if necessary, 
it may be assumed that no $W_k$ contains the origin. 

Now assume that $A$ has the Wiman-Valiron property as in Definition \ref{wvprop} and take $r \in F_1$ with $r$ large, and
the corresponding $z_r$. 
Let $Z(z)$ be defined on $D(z_r, 3)$ as in (\ref{Nlim5}). Then (\ref{b7}), (\ref{b8}) and
Lemma \ref{slem1} show that there exists 
$d > 0$, depending on $f$ but not on $k$ or $r$, such that 
\begin{eqnarray}
 \label{x1}
\chi( f(z), a_{j(k)}) &\geq& d \exp\left( - 20 |z| \sqrt{ |S_f(z)| } \right)  
\geq d \exp\left( - 30 |z| \sqrt{ |A(z)| } \right) \nonumber \\
&\geq& \exp\left( - 20 N |Z(z)|  \right) 
\geq \exp\left( - 30 N Z_r  \left| \frac{z}{r} \right|^{(N+2)/2}  \right), \quad N=N(r),
\end{eqnarray}
for all $z \in D(z_r, 3) \cap W_k$.


Next, take one of the islands $H_q$ in $D(z_r, 3)$ which by Lemma \ref{wvlema}
are mapped univalently by $Z(z)$ onto
$J_1$ as in (\ref{b8}), and write
\begin{equation}
 \label{x2}
f = \frac{C_1 u_1 + C_2 u_2}{D_1 u_1 + D_2 u_2} 
\end{equation}
for $z \in H_q$, in which $u_1$, $u_2$ are as in (\ref{b11})
and the $C_j$ and $D_j$ are complex numbers with $C_1D_2 - C_2 D_1 \neq 0$.
Since these coefficients may depend on $r$ and even on $q$, a number of lemmas are required in order to determine the behaviour of
$f(z)$ on certain arcs in $H_q$ which are pre-images of radial segments in $J_1$. 

\begin{lem}
 \label{segmentlem}
Let $K$ be a radial segment in $J_1$, given for some $a \in \R$ by 
\begin{equation*}
 \label{Kdefa}
K = \left\{ Z : \, a \leq \log \left| \frac{Z}{Z_r} \right| \leq a+1, \, \arg Z = \pm \frac{\pi}{8} \right\} \subseteq J_1.
\end{equation*}
Then on the pre-image $L= Z^{-1}(K)$ of $K$ in $H_q$ the $u_j$ satisfy 
\begin{equation}
 \label{x3a}
\pm \log \left| \frac{u_2(z)}{u_1(z)} \right| \geq  e^a  Z_r  \sin \frac{\pi}{8} .
\end{equation}
Moreover, if $r$ is large enough then  $L$ has  logarithmic diameter at least $1/N$. 
\end{lem}
\textit{Proof.} By (\ref{b11})  the $u_j$ and $U_j$ satisfy, for $z \in L$ and $Z \in K$,  
\begin{equation}
 \label{ujasymp}
\pm \log \left| \frac{u_2(z)}{u_1(z)} \right| =
\pm \log \left| \frac{U_2(Z)}{U_1(Z)} \right| = \pm {\rm Re} \, (2 i Z) + o(1) = 2 |Z| \sin \frac{\pi}{8} + o(1) ,
\end{equation}
and $|Z| \geq e^a Z_r$ 
on $K$. 
Furthermore, by (\ref{b7}), there exist points $h_1, h_2 \in L$ with
$$
\log \left| \frac{Z(h_2)}{Z(h_1)} \right| = 1, \quad 
\log \left| \frac{h_2}{h_1} \right| = \frac2{N+2} (1 + o(1) ) \geq \frac1{N} .
$$
\hfill$\Box$
\vspace{.1in}

\begin{lem}
 \label{lemx1}
There exists $j \in \{1, 2\}$ for which
$$\max \{ |C_1|, |C_2|, |D_1| , |D_2| \} = |C_j|.$$
\end{lem}
\textit{Proof.} Assume that this assertion is false: then  (\ref{x2}) may be written in the form
\begin{equation}
 \label{x4}
f = \frac{\alpha u + \beta v}{\gamma u + v} = \beta + \frac{(\alpha - \beta \gamma ) u/v}{\gamma u/v + 1} ,
\quad \{ u, v \} = \{ u_1, u_2 \}, 
\quad \max \{ |\alpha| , |\beta| , |\gamma| \} \leq 1.
\end{equation}
By (\ref{b8}) and Lemma \ref{segmentlem} there exists a radial segment $K^*$ in  $J_1$ given by 
\begin{equation}
 \label{Kdefb}
K^* = \left\{ Z : \, 0 \leq \log \left| \frac{Z}{Z_r} \right| \leq 1, \, \arg Z = \pm \frac{\pi}{8} \right\} ,
\end{equation}
on whose pre-image $L^*$ in $H_q$ the solutions satisfy  
\begin{equation}
 \label{x3b}
\log \left| \frac{u(z)}{v(z)} \right| \leq -  Z_r \sin \frac{\pi}{8}  ,
\end{equation}
and $L^*$ has logarithmic diameter at least $1/N$. 
But  (\ref{limZr}) and (\ref{x4}) now imply that
$$
\chi( f(z), \beta ) \leq |f(z)- \beta| \leq 4 \exp\left( -  Z_r \sin \frac{\pi}{8} \right) = o( N^{-2} ) 
$$
on $L^*$. Applying Lemma \ref{component} with $800 \pi \rho = 1/N$ shows that $f(z) - a_j$ is small on $L^*$, for some~$j$,
and this is a contradiction since $|\beta| \leq 1 < 2 \leq |a_j|$. 
\hfill$\Box$
\vspace{.1in}

It may henceforth be assumed that 
\begin{equation}
 \label{x5}
f = \frac{\alpha u + v}{\beta u + \gamma v } =   \frac{\alpha u/v + 1}{\beta u/v + \gamma },
\quad \max \{ |\alpha| , |\beta| , |\gamma| \} \leq 1, \quad \{ u, v \} = \{ u_1, u_2 \}. 
\end{equation}
\begin{lem}
 \label{lemx2}
The coefficient $\gamma$ satisfies $|\gamma| \geq N^{-3}$. 
\end{lem}
\textit{Proof.} Assume that this is not the case and take a radial segment $K^*$ defined as in (\ref{Kdefb}) such that
(\ref{x3b}) holds on the pre-image $L^*$ in $H_q$. 
Here $u$ and $v$ may not be the same choices of $u_1$ and $u_2$ as
in the previous lemma, but this does not affect the existence of $K^*$. 
Then on $L^*$ the representation (\ref{x5}) gives 
$$
\alpha u(z)/v(z) + 1 \sim 1, \quad |\beta u(z)/v(z) + \gamma | \leq 2 N^{-3}, 
\quad \chi(f(z), \infty) = o( N^{-2} ).$$ 
Applying Lemma \ref{component} again with $800 \pi \rho = 1/N$ shows that $f(z) - a_j$ is small on $H^*$, for some~$j$,
and this is a contradiction since  the $a_j$ are all finite. 
\hfill$\Box$
\vspace{.1in}

In view of Lemma \ref{lemx2}, the representation (\ref{x5}) can now be written in the form 
\begin{equation}
 \label{x6}
f = \frac1\gamma + \frac{\alpha-\beta/\gamma}{\beta + \gamma v/u}, \quad \left| \frac1\gamma \right| \leq N^3  .
\end{equation}
\begin{lem}
 \label{lemx3}
The coefficients in (\ref{x6}) satisfy
\begin{equation}
 \label{x7}
2 N^3 \geq | \alpha-\beta/\gamma | \geq \exp( - 400 N Z_r  ) .
\end{equation}
\end{lem}
\textit{Proof.} In (\ref{x7}) the first inequality follows at once from (\ref{x5}) and (\ref{x6}), so assume that the second inequality is false.
Again take a radial segment $K^*$ defined as in (\ref{Kdefb}) such that
(\ref{x3b}) holds on the pre-image $L^*$ in $H_q$, and recall that $L^*$ has logarithmic diameter at least $1/N$. 
By (\ref{limZr}) this yields, for $z \in L^*$, since $| \beta | \leq 1$ but $|\gamma| \geq N^{-3}$,
$$
| \beta + \gamma v(z)/u(z)| \geq 1, \quad \chi( f(z), 1/\gamma ) \leq |f(z) - 1/\gamma| < \exp( - 400  N Z_r ) .
$$
Applying Lemma \ref{component}, with $\rho = \exp( - 200  N Z_r )  =o(1/N)$, gives 
$$
\chi(  f(z) , a_j ) \leq \exp( -  100  N Z_r ) 
$$
on $L^*$, for some $j$, and Lemma \ref{lemp2} shows that the component of the set $\{ z \in \C : |f(z)-a_j| < \tau \}$
which contains $L^*$ cannot also contain a zero of $f-a_j$. Thus
$L^*$ must lie in one of the logarithmic tracts $W_k$, with corresponding 
asymptotic value $a_{j(k)} = a_j$. Now 
(\ref{x1}) and
(\ref{Kdefb}) imply that there exists $z \in L^*$ with
$$
|Z(z)| = Z_r, 
\quad \chi( f(z), a_j) \geq  \exp\left( - 20 N Z_r    \right)  ,
$$
and this is a contradiction.
\hfill$\Box$
\vspace{.1in}

\begin{lem}
 \label{lemx4}
The coefficients  in (\ref{x6}) satisfy
\begin{equation}
 \label{x8}
\left| \frac{\alpha}{\beta} \right|   \leq \left| \frac{1}{\beta} \right| \leq \exp( 800 N Z_r  ) .
\end{equation}
\end{lem}
\textit{Proof.} The first inequality is obvious since $| \alpha | \leq 1$, so assume that the second fails.
Take  a radial segment $K^{**}$ in  $J_1$ given by 
\begin{equation*}
 \label{Kdefc}
K^{**} = \left\{ Z : \, 3 \log N  \leq \log \left| \frac{Z}{Z_r} \right| \leq 1 + 3 \log N , \, \arg Z = \pm \frac{\pi}{8} \right\} ,
\end{equation*}
on whose pre-image $L^{**}$ in $H_q$ the solutions satisfy  
\begin{equation*}
 \label{x3bc}
 \log \left| \frac{u(z)}{v(z)} \right| \geq N^3 Z_r \sin \frac{\pi}{8} , \quad 
N^2 Z_r = o \left( \log \left| \frac{u(z)}{v(z)} \right| \right)  .
\end{equation*}
This is possible by (\ref{Nlim}), (\ref{b8}) and (\ref{x3a}), and $L^{**}$ has logarithmic diameter at least $1/N$. 
For $z \in L^{**}$ this yields, in view of (\ref{x6}) and (\ref{x7}),  
$$
| \beta + \gamma v(z)/u(z)| \leq 2 \exp( - 800 N Z_r  ), \quad 
\left| \frac{\alpha-\beta/\gamma}{\beta + \gamma v(z)/u(z)} \right| \geq  \exp( 200 N Z_r  ), $$
and so $|f(z)| \geq   \exp( 100 N Z_r  ) $, which gives a contradiction as in the proof of Lemma \ref{lemx2}. 
\hfill$\Box$
\vspace{.1in}

Since $\beta \neq 0$, the representation (\ref{x5}) for $f$ can also be written in the form
\begin{equation}
 \label{x8a}
f = \frac{\alpha}{\beta}  + \frac{1 - \alpha \gamma /\beta}{\beta u/v  + \gamma }  ,
\end{equation}
and it follows at once from (\ref{x5}) and (\ref{x7}) that
\begin{equation}
 \label{x9}
\left| \frac{\alpha }{\beta} - \frac1{\gamma} \right| = \left| \frac{1 }{\beta}\right| \cdot
\left| \alpha  - \frac{\beta}{\gamma} \right| \geq \exp( - 400 N Z_r ). 
\end{equation}
In particular, $\alpha /\beta \neq 1/\gamma$. Since $r$ is large the next lemma is an immediate consequence of (\ref{x1}). 

\begin{lem}
 \label{lemx6a}
Let 
\begin{equation}
 \label{Rrdef}
R_r = r \exp \left( \frac{20 \log N}{N+2} \right) . 
\end{equation}
For each  tract $W_k$ and its associated asymptotic value $a_j = a_{j(k)}$, 
the function $f$ satisfies  
\begin{equation}
 \label{x88a}
|f(z)-a_j| \geq \chi(f(z), a_j) \geq   \exp \left( - N^{12} Z_r  \right) \quad \hbox{for $z \in D(z_r, 3) \cap W_k$ with $|z| \leq R_r$}.
\end{equation}
\end{lem}
\hfill$\Box$
\vspace{.1in}
\begin{lem}
 \label{lemx5}
There exist disjoint radial segments $K^+, K^-$ in $J_1$, each given by 
\begin{equation}
 \label{x10}
K^{\pm} = \left\{ Z : \, 5 \log N  \leq \log \left| \frac{Z}{Z_r} \right| \leq \frac12 \phi(r) , \, | \arg Z | =  \frac{\pi}{8} \right\} ,
\end{equation}
on whose  pre-images $L^{\pm}$ in $H_q$ the function $f$ satisfies 
\begin{equation}
\label{x11}
\chi(f(z), 1/\gamma) \leq |f(z) - 1/\gamma | \leq \exp \left( - |Z| \sin \frac{\pi}8 \right) \leq \exp \left( - N^3 Z_r \right) 
\quad \hbox{for $z \in L^+$,} 
\end{equation}
\begin{equation}
\label{x12}
\chi(f(z), \alpha/\beta ) \leq
|f(z) - \alpha /\beta  | \leq \exp \left( - |Z| \sin \frac{\pi}8 \right)  \leq \exp \left( - N^3 Z_r \right) 
\quad \hbox{for $z \in L^-$.} 
\end{equation}
Moreover, there exist distinct singular values $a_{j+} $ and $a_{j-}$ of $f^{-1}$, such that
$L^{\pm}$ lies in a logarithmic tract $W_{k \pm}$ with  asymptotic value $a_{j \pm}$. 
With a fixed determination of $\arg Z(z_r)$, each pre-image $L^{\pm}$ is a simple arc in $H_q$, on which
\begin{equation}
 \label{argz}
\arg \frac{z}{z_r} = \frac{2}{N+2} \left(\pm \frac\pi8 + 2 \pi P_q - \arg Z(z_r) \right)   + o \left( \frac1{N} \right) = o(1), \quad P_q \in \Z, 
\end{equation}
and distinct $H_q$ give rise to distinct integers $P_q$. Moreover, $L^{\pm}$ meets the circle $S(0, R_r)$ non-tangentially 
and contains a simple arc $M^{\pm}$ which lies in $|z| \geq R_r$ and joins the unique point 
in $L^{\pm} \cap S(0, R_r)$ to a point 
 $\zeta_{\pm}$ which satisfies  $| \zeta_{\pm} | > R_r$ and
\begin{equation}
 \label{zetadef}
|f(z) - a_{j\pm} | >  | f(\zeta_{\pm}) - a_{j\pm} | = \exp \left( - N^{100} Z_r \right) 
\quad \hbox{ for all $z \in M^{\pm} \setminus \{ \zeta_{\pm} \}$.}
\end{equation}
\end{lem}
\textit{Proof.} The segment $K^-$ is chosen so that, using (\ref{ujasymp}) and (\ref{x10}), 
$$
\log \left| \frac{u(z)}{v(z)} \right| = 2 |Z| \sin \frac{\pi}8 + o(1) \geq N^4 Z_r .
$$
On combination with (\ref{x5}), 
(\ref{x8}) and (\ref{x8a}) this proves (\ref{x12}), and (\ref{x11}) is derived from (\ref{x6}) in the same way.
Since $L^{\pm}$ has logarithmic diameter at least $1/N$, by
Lemma \ref{segmentlem}, applying Lemma \ref{component} gives singular values $a_{j+}$, $a_{j-}$ of $f^{-1}$ with 
$$
\chi(f(z), a_{j+} ) \leq \exp \left( - N^2 Z_r \right)  \quad \hbox{for $z \in L^+$}, \quad 
\chi(f(z), a_{j-} ) \leq \exp \left( - N^2 Z_r \right)  \quad \hbox{for $z \in L^-$}, 
$$
and so
$$
\max \{ |\alpha/\beta | , |1 /\gamma | \} \leq \max \{ |a_j| \} + 1.
$$
It now follows using (\ref{x9}) that $a_{j+} \neq a_{j-}$, and the inclusion $L^{\pm} \subseteq W_{k \pm}$ is a consequence of Lemmas \ref{lemp2}
and \ref{segmentlem}. 

To prove the remaining assertions, write $Z = Z_r e^{T+iB}$ on $K^{\pm}$, where $|B | =  \pi /8$
and $T$ is real. The formulas (\ref{b7}) and (\ref{b8}) imply that there exists $E_q \in \C$ such that
$$
\frac{N+2}2  \log \frac{z}{z_r} =  \log \frac {Z}{Z_r} + E_q  + o(1)  \quad
\hbox{for $z \in H_q$ and $Z = Z(z) \in J_1$,}
$$
in which the logarithms are chosen to vanish at $z_r$ and $Z_r$ respectively. Taking exponentials shows that 
$E_q = i( 2 \pi P_q - \arg Z(z_r) ) + o(1)$ for some integer $P_q$, and there exists $z \in H_q$ with
$$
\frac{N+2}2  \log \frac{z}{z_r} =  i( 2 \pi P_q - \arg Z(z_r) ) ;
$$
this follows from (\ref{b7}), Koebe's quarter theorem and the fact that
$\log z/z_r $ is a univalent function of
$\log Z$ for $Z $ in  $ J_1$ by Lemma \ref{wvlema}. 
Hence distinct $H_q$ must give rise to distinct integers $P_q$, and 
(\ref{argz}) follows for $z \in L^{\pm}$ and $Z = Z_r e^{T+iB}$ on $K^{\pm}$; here
$\arg (z/z_r) = o(1)$ by (\ref{Nlim}) and (\ref{Nlim1}).

The arcs $L^{\pm}$ meet $S(0, R_r)$, by (\ref{Nlim}), (\ref{b7}), (\ref{b8}), (\ref{Rrdef}) and the fact that
$T$ varies between $5 \log N$ and $\phi(r)/2$,
and they do so non-tangentially, because (\ref{Nlim5}) and (\ref{b7}) give
\begin{equation}
 \label{dzdT}
\frac{dZ}{dz} = A(z)^{1/2} \sim \frac{(N+2)Z}{2z} , \quad
\frac{dz}{dT} = \frac{dz}{dZ} \, \frac{dZ}{dT} = \frac{dz}{dZ} \, Z \sim \frac{2z}{N+2} .
\end{equation}
Observe next that the pre-image $Q^{\pm}$ in $H_q$ of the part of $K^{\pm}$ with $\phi(r)/2 - 1 \leq T \leq \phi(r)/2$ 
has logarithmic diameter at least $1/N$, by Lemma \ref{segmentlem}, and 
$$
\chi( f(z), \delta ) \leq  \exp \left( - |Z| \sin \frac{\pi}8 \right) 
\leq  \exp \left( - N^{300} Z_r \right)  \quad \hbox{on $Q^{\pm}$, where} \quad \delta \in \{ \alpha/\beta, 1/\gamma \}, 
$$
by (\ref{Nlim}), (\ref{x11}) and (\ref{x12}). Applying Lemma~\ref{component} and (\ref{x88a}) shows that
\begin{equation}
 \label{Qest}
|f(z)-a_{j\pm} |\leq \exp \left( - N^{100} Z_r \right) \quad \hbox{and} \quad  |z| > R_r  
\quad \hbox{for} \quad 
z \in Q^{\pm} \subseteq L^{\pm} \subseteq W_{k \pm}. 
\end{equation}

The arc $M^{\pm}$ and the point $\zeta_{\pm} $ are now determined as follows. Since (\ref{dzdT}) shows that 
$$
\frac1z \frac{dz}{dT} = \frac{ d \log z}{dT}
$$
has positive real part, it follows that $|z|$ increases with $T$ as $Z$ describes $K^{\pm}$. 
Start at the unique point in $L^{\pm} \cap S(0, R_r)$ and follow $L^{\pm}$ in the direction of increasing $T$, 
and let $\zeta_{\pm}$ be the  first point encountered
which satisfies $|f(z) - a_{j\pm} |  = \exp \left( - N^{100} Z_r \right)$, 
such a point existing because of~(\ref{Qest}).
\hfill$\Box$
\vspace{.1in}

There are finitely many
logarithmic tracts $W_k$, but arbitrarily many islands $H_q \subseteq D(z_r, 3)$, each with corresponding 
tracts 
$W_{k+} \neq W_{k-}$, which makes it possible to choose tracts and arcs as follows. Take logarithmic tracts
$U_a$ and $U_b$, with distinct corresponding asymptotic values $a$ and $b$ respectively, and 
disjoint arcs $M_{p,a}$  contained in $U_a \cap D(z_r, 3)$, as well as
disjoint arcs $M_{p,b}$ contained in $U_b \cap D(z_r, 3)$, in each case for $p=1, 2$.
Each of these arcs is one of the $M^{\pm}$ as in Lemma \ref{lemx5}. 
Here $M_{p,a}$ joins a point $x_{p,a}$ on $S(0, R_r)$ to a point $y_{p,a}$ satisfying
$$| y_{p,a} | > R_r , 
\quad | f(y_{p,a}) - a | = \exp \left( - N^{100} Z_r \right),
$$
and the $M_{p,b}$ join corresponding points $x_{p,b}$ and  $y_{p,b}$. In view of 
(\ref{argz}) and the fact that the number of islands $H_q$ may be chosen arbitrarily large, it may be assumed that 
as $z$ describes $D(z_r, 3) \cap S(0, R_r)$ with $\arg z$ increasing, $z$ meets these arcs in the order
$M_{1,a}$, $M_{1,b}$, $M_{2,a}$, $M_{2,b}$. 

Now $y_{1,a}$ and $y_{2,a}$ both lie in $U_a$, and $\log 1/(f-a)$ maps $U_a$ univalently onto the half-plane
${\rm Re} \, w > - \log \tau $. Thus there exists a simple path $X_1 \subseteq U_a$ joining 
$y_{1,a}$ to $y_{2,a}$, on which 
$$
| f(z) - a | = \exp \left( - N^{100} Z_r \right) ,
$$
and $X_1$ does not meet $D(z_r, 3) \cap S(0, R_r)$, because of (\ref{x88a}). Recalling 
(\ref{zetadef}), it is now possible to form a Jordan curve $X$ consisting of
$X_1$, $M_{1,a}$, $M_{2,a}$ and an arc $X_2$ of $S(0, R_r) \cap D(z_r, 3)$ joining $x_{1,a}$ to $x_{2,a}$, in such a way that
 $X \cap S(0, R_r) \cap D(z_r, 3)
= X_2$. Now join $x_{1, b}$ to $x_{2,b}$ by a simple path $Y_1$ which, apart from its end-points, lies in $B(0, R_r) \cap D(z_r, 3)$; this 
may be done so that $Y_1$  
meets $S(0, R_r)$ non-tangentially at $x_{1, b}$ and $x_{2,b}$, and  meets $X$ only at $x_{1,b}$.  
This gives a simple path $Y $ from $y_{1,b}$ to $y_{2,b}$,
consisting of $Y_1$, $M_{1,b}$ and $M_{2,b}$, which intersects $X$ only once, at $x_{1,b}$, and does so non-tangentially.
Hence one of $y_{1,b}$ and $y_{2,b}$, say $y_b$,  lies in the bounded component of $\C \setminus X$. But $y_b$ lies in $U_b$, 
and can be joined to infinity by a path in $U_b$ on which  
$$
| f(z) - b | \leq \exp \left( - N^{100} Z_r \right) .
$$
This path cannot meet $X$, because of (\ref{x88a}), which gives a contradiction and completes the proof
of Theorem \ref{thm2}.
\hfill$\Box$
\vspace{.1in}

\noindent
\textit{Remark.} 
Examples such as $\exp( \exp(z))$ make it clear that hypothesis (ii) is not redundant in Theorems \ref{thm1} and \ref{thm2}.
It is reasonable to ask whether the assumption  on the multiplicities could be eliminated in hypothesis (i), but
difficulties would appear to arise in the application of Lemma \ref{lemp2}. Here the example $f_n(z) = (1+z/n)^n $, for which
$S_{f_n}$ converges locally uniformly to $-1/2 = S_{\exp}$, suggests that in general it is not easy to deduce information about
multiplicities from the growth of the Schwarzian.

\section{The case of Borel exceptional poles}

Following the terminology of Hayman \cite{delta,Hay8}, define an $\varepsilon$-set to be a countable
union of discs 
\begin{equation}
E = \bigcup_{j=1}^\infty B(b_j, r_j) \quad \hbox{ such that} \quad
\lim_{j\to \infty} |b_j| = \infty \quad \hbox{and} \quad
\sum_{j=1}^\infty \frac{r_j}{|b_j|} < \infty .
\label{car1}
\end{equation}
Note that if $E$ is an $\varepsilon$-set then the set of $r \geq 1$ 
for which the circle $S(0, r)$ meets $E$ has finite logarithmic 
measure. The next two lemmas are refinements of estimates from \cite{delta}. 

\begin{lem}\label{lem1}
Let $\alpha > 1$ and $\lambda > 0$. 
Let $b_1, b_2, \ldots $ be complex numbers with
$|b_k| \leq |b_{k+1}| $ and $\lim_{k\to\infty} |b_k| = \infty$, as well as
$\sum_{b_k \neq 0} |b_k|^{-\lambda} < \infty$.
For $r > 0$ let $n(r)$ be the number of $b_k$, 
taking account of repetition, with $|b_k| \leq r$.  
Then there exist a positive constant
$d_\alpha$, depending only on $\alpha$, and
an $\varepsilon$-set $E = E_\alpha$ such that for large $z$ with
$z \not \in E$ and $|z| = r$,
\begin{equation}
\sum_{|b_k| < \alpha r } \frac1{|z - b_k|} <
d_\alpha \frac{n( \alpha^2 r) }r ( \log r)^\alpha \log n( \alpha^2 r) 
\label{car2}
\end{equation}
and 
$|z - b_k| \geq 4 |b_k|^{1-\lambda} $ for all $k$.
\end{lem}
{\em Proof.} 
By \cite[formulas (5.6) and (5.8)]{Gun2}, the estimate (\ref{car2}) holds outside an exceptional set $E$
satisfying (\ref{car1}).
Moreover, for large $k_0$  the set
$E' = \bigcup_{k \geq k_0} B(b_k, 4 |b_k|^{1-\lambda}) $ is an $\varepsilon$-set,
and it suffices to replace $E$ by $E \cup E'$.
\hfill$\Box$
\vspace{.1in}

\begin{lem}\label{lem3}
Let $G$ be a transcendental meromorphic function in the plane 
with order of growth less than~$\lambda < \infty $. Then there exists an
$\varepsilon$-set $E$ such that, as $z \to \infty$ in $\C \setminus E$,
\begin{equation}
G( \zeta ) \sim G(z)  \quad \hbox{and} \quad | \log |G(\zeta)| | = o\left( |z|^\lambda  \right) \quad
\hbox{ for } \quad | \zeta - z| \leq |z|^{1 - \lambda} .
\label{car5}
\end{equation}
\end{lem}
{\em Proof.} Since $G$ has order less than $\lambda$ the sequence $(b_k)$ of
zeros and poles of $G$,  repeated according to multiplicity, 
satisfies $\sum_{b_k \neq 0} |b_k|^{-\lambda} < \infty$. Applying Lemma 
\ref{lem1}  with $\alpha = 4$ 
gives an $\varepsilon$-set $E$ such that
for large $z \not \in E$ the estimate     
(\ref{car2}) holds,
as well as $|z - b_k| \geq 4 |b_k|^{1-\lambda}$ for all~$k$. 
Let $z$ be large, 
not in $E$, set $r = |z|$, and let $|u - z | \leq |z|^{1 - \lambda}$. Then, for all $k$, 
\begin{equation}
r/2 \leq | u | \leq 2 r 
\quad \hbox{and} \quad 
| u - b_k | \geq | z - b_k | - |z|^{1 - \lambda} 
\geq \frac{ | z - b_k | }2  > 0  .
\label{car6}
\end{equation}
Now (\ref{car2}), (\ref{car6})
and the standard estimate for $|G'/G|$ (see \cite[p.65]{Jank}) give $d > 0$ such that
\begin{eqnarray*}
\left| \frac{G'(u)}
{G(u)} \right| 
&\leq& d \frac{T(2|u|, G)}{|u|} + \sum_{|b_k| < 2|u|  } \frac2{|u - b_k|} 
\leq 2d \frac{T(4 r, G)}r + \sum_{|b_k| < 4r } \frac4{|z - b_k|} \\
&\leq& 2d \frac{T(4 r, G)}r + 4 d_4 
\frac{n(16r)}{r} (\log r)^4 \log n(16r)
= o\left(r^{\lambda - 1}\right) ,
\end{eqnarray*}
where $n(r) = n(r, G) + n(r, 1/G)$.
The  first assertion of (\ref{car5}) now follows by integration of $G'/G$, and the second is a routine consequence of the Poisson-Jensen formula. 
\hfill$\Box$
\vspace{.1in}

\begin{lem}\label{lemwv1}
Let $A$ be a transcendental meromorphic function in the plane for which $\infty$ is a Borel exceptional value. 
Then $A$ has the Wiman-Valiron property.
\end{lem}
\textit{Proof.} 
There is nothing to prove if $A$ has finitely many poles, so assume that this is not the case, and write 
\begin{equation}
 \label{n1}
A = \frac{H}{G},
\end{equation}
in which $H$ and $G$ are transcendental entire functions with $G \not \equiv 0$ and
$\rho(G) < \rho(A) = \rho(H)$. Choose $\lambda$ and
$\gamma$ with
\begin{equation}
 \label{n2}
\frac12 < \gamma < 1, \quad \rho( G) < \lambda < \gamma \rho (H) <   \rho (H) . 
\end{equation}
By standard facts from the Wiman-Valiron theory \cite{Hay5} there exists
a set $F_2 \subseteq [1, \infty )$ with the following properties. First,
$[1, \infty) \setminus F_2$ has finite logarithmic measure. Second, the central index $N(r)$ of $H$ satisfies
\begin{equation}
 \label{n3}
N(r)^\gamma \leq  \log M(r, H) 
\end{equation}
for $r \in F_2$. Finally, if  $|z_r| = r \in F_2$ and $|H(z_r)| \sim M(r, H)$ then
$$
H(z) \sim \left( \frac{z}{z_r} \right)^{N(r)} H(z_r) \quad \hbox{for} \quad
 z \in E_r = 
\left\{  z_r e^\tau : \, 
| {\rm Re} \, \tau  | \leq 8 N(r)^{-\gamma} \, , \,  
| {\rm Im } \, \tau | \leq 8  N(r)^{-\gamma}  \, \right\}.
$$

It may be assumed in addition that for $r \in F_2$ the circle $S(0, r)$ does not meet the $\varepsilon$-set $E$ of Lemma \ref{lem3}. 
Since
$N(r)$ has the same order of growth $\rho(H)$ as $\log M(r, H)$, there exists  by (\ref{n2}) a set $F_1 \subseteq F_2$,
of infinite logarithmic measure, such that 
\begin{equation}
\label{n8}
r^\lambda = o( N(r)^{\gamma} ) \quad \hbox{for} \quad r \in F_1. 
\end{equation} 
In particular, if $z \in E_r$  then, by  (\ref{n8}),  
$$|z-z_r| = O( r N(r)^{ - \gamma } ) = 
o( r^{1-\lambda}) $$
and so (\ref{car5}), (\ref{n1}), (\ref{n3}) and (\ref{n8}) yield $G(z) \sim G(z_r)$
and 
$$\log |G(z)| = o( r^\lambda ) = o ( N(r)^{  \gamma } ) = o( \log M(r, H)), \quad \log |A(z_r)| \sim \log M(r, H).$$ 
Thus setting $\phi(r) =  N(r)^{1-\gamma} $ gives
$$
\frac{\phi(r)}{ \log^+ |A(z_r)| } \leq  2 N(r)^{1-2 \gamma} \to 0
$$
as $r \to \infty$ in $F_1$, using (\ref{n2}). It follows that $F_1$, $z_r$ and $A$ satisfy Definition \ref{wvprop}. 
\hfill$\Box$
\vspace{.1in}

\noindent
School of Mathematical Sciences, University of Nottingham, NG7 2RD.\\
jkl@maths.nott.ac.uk


{\footnotesize
\begin{thebibliography}{99}
\bibitem{BE}
W. Bergweiler and A. Eremenko, On the singularities of the inverse to
a meromorphic function of finite order,
\textit{Rev. Mat. Iberoamericana } 11 (1995), 355-373.
\bibitem{delta}
W. Bergweiler and J.K. Langley,
Zeros of differences of meromorphic functions,
\textit{Math. Proc. Camb. Phil. Soc.} 142 (2007), 133-147.
\bibitem{BRS}
W. Bergweiler, P.J. Rippon and G.M. Stallard,
Dynamics of meromorphic functions with direct or logarithmic
singularities, , \textit{Proc. London Math. Soc.} 97 (2008), 368-400.
\bibitem{Dra1}
D. Drasin,
Proof of a conjecture of F. Nevanlinna concerning functions
which have deficiency sum two, \textit{Acta. Math.} 158 (1987), 1-94.
\bibitem{elf}
G. Elfving, \"Uber eine Klasse von Riemannschen Fl\"achen und ihre Uniformisierung, \textit{Acta Soc. Sci. Fenn.} 2 (1934) 1-60.
\bibitem{Er2}
A. Eremenko, Meromorphic functions with small ramification,
\textit{Indiana Univ. Math. J.} 42 (1994), 1193-1218.
\bibitem{Erem}
A. Eremenko,
Geometric theory of meromorphic functions, 
\textit{In the Tradition of Ahlfors-Bers, III},
\textit{Contemp. Math. } 355, (Amer. Math. Soc. Providence, 2004) 221-230. 
\bibitem{EL}
A. Eremenko and M.Yu. Lyubich, Dynamical properties of some classes
of entire functions, \textit{Ann. Inst. Fourier Grenoble} 42 (1992), 989-1020.
\bibitem{GO} A.A. Gol'dberg and I. V. Ostrovskii,
\textit{Distribution of values of meromorphic functions}, (Nauka, Moscow, 1970 (Russian)).
English transl., \textit{Translations of Mathematical Monographs} 236, (Amer. Math. Soc. Providence 2008).
\bibitem{Gun2}
G. Gundersen, Estimates for the logarithmic derivative
of a meromorphic function, plus similar estimates, \textit{J.~London Math. Soc.}
(2) 37 (1988), 88-104.
\bibitem{Hay2} W.K. Hayman, \textit{Meromorphic functions},  (Oxford at the Clarendon Press, 1964).
\bibitem{Hay5}
W.K. Hayman, The local growth of power series: a survey of the
Wiman-Valiron method, \textit{Canad. Math. Bull. } 17 (1974), 317-358.
\bibitem{Hay8}
W.K. Hayman, Slowly growing integral and subharmonic functions,
\textit{Comment. Math. Helv.} 34 (1960), 75-84.
\bibitem{Hay9}
W.K. Hayman, \textit{Multivalent functions}, 2nd edition, \textit{Cambridge Tracts in
Mathematics} 110, (Cambridge University Press, Cambridge, 1994).
\bibitem{Hil1}
E. Hille, \textit{Lectures on ordinary differential equations}, (Addison-Wesley,
Reading, Mass., 1969).
\bibitem{Hil2}
E. Hille, \textit{Ordinary differential equations in the complex domain}, (Wiley,
New York, 1976).
\bibitem{Iversen}
F. Iversen, \textit{Recherches sur les fonctions inverses des fonctions m\'eromorphes}, Thesis, University of Helsinki, 
(Imprimerie de la Soci\'et\'e de Litt\'erature Finnoise, Helsinki, 1914).
\bibitem{Jank}
G. Jank and L. Volkmann, \textit{Einf\"uhrung in die Theorie der
ganzen und meromorphen Funktionen mit Anwendungen auf
Differentialgleichungen}, (Birkh\"auser, Basel, 1985).
\bibitem{La5}
J.K. Langley, Proof of a conjecture of Hayman concerning $f$ and
$f'' $, \textit{J. London Math. Soc. } (2) 48 (1993), 500-514.
\bibitem{Nehari}
Z. Nehari, 
The Schwarzian derivative and schlicht functions.
\textit{Bulletin Amer. Math. Soc.} 55 (1949), 545-551. 
\bibitem{Nehari2}
Z. Nehari, \textit{Conformal mapping}, (Dover, New York, 1975). 
\bibitem{FNev}
F. Nevanlinna,\"Uber eine Klasse meromorpher Funktionen, 
\textit{Septi\`eme Congr\`es Math. Scand.} Jbuch. 56 (Oslo, 1929), 81-83.
\bibitem{Nev2}
R. Nevanlinna, \"Uber Riemannsche Fl\"achen mit endlich vielen Windungspunkten,
\textit{ Acta Math.} 58 (1932), 295-373.
\bibitem{Nev}
R. Nevanlinna, \textit{Eindeutige analytische Funktionen,
2. Auflage } (Springer, Berlin, 1953).



\end{thebibliography}
}
\end{document}